\documentclass[11pt,a4paper]{article}

\usepackage{amsmath,amssymb}
\newtheorem{theorem}{Theorem}[section]

\newtheorem{definition}{Definition}[section]
\newtheorem{remark}{Remark}[section]
\numberwithin{equation}{section}
\newtheorem{example}{Example}[section]
\newenvironment{proof}[1][Proof]{\noindent\textbf{#1.} }{\hfill {$\Box$}}
\numberwithin{equation}{section}

\begin{document}

\title{\textsc{Pointwise Trichotomy for Skew-Evolution Semiflows on Banach Spaces }}
\author{\textsc{Codru\c{t}a Stoica}}
\date{}
\maketitle

{\footnotesize \noindent \textbf{Abstract.} The paper introduces
the notion of skew-evolution semiflows and presents the concept of
pointwise trichotomy in the case of skew-evolution semiflows on a
Banach space. The connection with the classic notion of trichotomy
presented in \cite{MeSt_IEOT}, for evolution operators, is also
emphasized, as well as some characterizations. The approach of the
theory is from uniform point of view. The study can also be
extended to systems with control whose state evolution can be
described by skew-evolution semiflows.}

{\footnotesize \vspace{3mm} }

{\footnotesize \noindent \textit{Mathematics Subject
Classification:} 34D09}

{\footnotesize \vspace{2mm} }

{\footnotesize \noindent \textit{Keywords:} Skew-evolution
semiflow, exponential trichotomy, pointwise exponential
trichotomy}

\section{Introduction}

It is in \cite{ChLe_JDE} that S.N. Chow and H. Leiva introduced
and characterized the notion of pointwise discrete dichotomy for
discrete skew-product flows. The existence of this asymptotic
property is related to the existence of a unique bounded solution
of  certain difference equation, involving the discrete cocycle.
They also established the connection between the latest concept
and the classical uniform exponential dichotomy for linear
skew-product flows.

In \cite{LaMoRa_JDE} Y. Latushkin, S. Montgomery-Smith and T.
Randolph studied the spectral mapping and the hyperbolicity,
related to the exponential dichotomy, for evolutionary semigroups
generated by a strongly continuous semicocycle over a locally
compact metric space acting on Banach fibers. The results are
applied in the study of pointwise dichotomy versus global
dichotomy. Later, in \cite{LaSc_JDE} is studied the exponential
dichotomy of an exponentially bounded, strongly continuous cocycle
on a locally compact metric space acting on a Banach space.
Dichotomy is characterized in terms of hyperbolicity of a family
of weighted shift operators. It is shown that exponential
dichotomy follows from pointwise discrete dichotomies with uniform
constants.

In \cite{MeSaSa_IEOT} M. Megan, A.L. Sasu and B. Sasu emphasize
the concept of pointwise exponential dichotomy for linear
skew-product flows, by necessary and sufficient conditions, as
well as the relation with global exponential dichotomy. It is
proved that the presented property for skew-product flows is
equivalent with the pointwise admissibility of the pair
$(C_{0}(\mathbb{R}, \mathcal{V}), C_{0}(\mathbb{R},
\mathcal{V}))$, $\mathcal{V}$ being a Banach space.

Systems with control are studied in \cite{MeSaSa_RMC}, their state
evolution being described by linear skew-product semiflows. It is
also presented the relationship between the concepts of exact
controllability and complete stabilizability of general control
systems.

The concept of trichotomy is considered as a natural
generalization of the classic concept of dichotomy. The basic idea
is to obtain, at any moment, a splitting of the space into three
subspaces: the stable subspace, the unstable subspace and the
neutral one. The notion of Sacker-Sell type trichotomy was
introduced in 1976 and studied in the case of linear differential
equations in \cite{SaSe_JDE}. For the first time a sufficient
condition for the existence of trichotomy in the case of
$\mathbb{R}^{n}$ was emphasized.

Later, in 1988, a stronger notion was introduced by S. Elaydi and
O. Hajek in \cite{ElHa_JMAA}, namely the exponential trichotomy
for linear and nonlinear differential systems by means of Liapunov
functions.

In \cite{MeSt_IEOT} we have presented the property of uniform
exponential trichotomy in the case of evolution operators in
Banach spaces. The results extend well-known theorems obtained in
the linear case for uniform exponential stability, but we
emphasize a unified treatment for uniform asymptotic behaviors
(exponential decay, exponential growth, exponential stability,
exponential instability, exponential dichotomy, exponential
trichotomy, as given in \cite{MeStBu_UVT}) in the setting of
nonlinear evolution operators.

Equivalent definitions and characterizations for the uniform
exponential trichotomy by means of two respectively four
compatible projector families are given in \cite{MeSt_OT}. Some
properties for a particular case of skew-evolution semiflows were
introduced in \cite{MeStBu_CJM}.

The theory of skew-evolution semiflows introduced by us is used to
study the asymptotic behavior of time-varying linear systems, the
question that arises being whether the connection between
stabilizability and controllability can be extended to systems
associated to skew-evolution semiflows.

In our present paper we extend the study of exponential trichotomy
in a more general setting fitting the new concept of
skew-evolution semiflows, from a global and a pointwise point of
view, some similarities and differences being emphasized. The
study of the concept of skew-evolution semiflow for the approach
of evolution equations by means of the evolution operators theory
was essential.

\section{Definitions and examples}

Let $(\mathcal{X},d)$ be a metric space, $\mathcal{V}$ a Banach
space, $\mathcal{B}(\mathcal{V})$ the space of all
$\mathcal{V}$-valued bounded operators defined on $\mathcal{V}$.
We denote $\mathcal{Y}=\mathcal{X}\times \mathcal{V}$ and we
consider the set
$\mathcal{T}=\left\{ (t,t_{0})\in \mathbb{R}%
_{+}^{2}:t\geq t_{0}\right\}$. The norm of vectors on
$\mathcal{V}$ and operators on $\mathcal{B}(\mathcal{V})$ is
denoted by $\left\Vert \cdot \right\Vert$. Let $I$ be the identity
operator on $\mathcal{V}$.

\begin{definition}\rm\label{d_semiflow}
A mapping $\psi:\mathcal{T}\times \mathcal{X}\rightarrow
\mathcal{X}$ is called \textit{evolution semiflow} on
$\mathcal{X}$ if it satisfies the following properties

$(es_{1})$ $\psi(t,t,x)=x, \ \forall (t,x)\in
\mathbb{R}_{+}\times \mathcal{X}$

$(es_{2})$ $\psi(t,s,\psi(s,t_{0},x))=\psi(t,t_{0},x), \
\forall (t,s),(s,t_{0})\in \mathcal{T}, \ \forall x\in
\mathcal{X}$.
\end{definition}

\begin{definition}\rm\label{d_cociclu_2}
A mapping $\Psi:\mathcal{T}\times \mathcal{X}\rightarrow
\mathcal{B}(\mathcal{V})$ that satisfies the following properties

$(ec_{1})$ $\Psi(t,t,x)=I, \ \forall t\geq0,\ \forall x\in
\mathcal{X}$

$(ec_{2})$ $\Psi(t,t_{0},x)=\Psi(t,s,\psi(s,t_{0},x))\Psi(s,t_{0},x),
\forall (t,s), (s,t_{0})\in \mathcal{T}, \forall x\in
\mathcal{X}$

\noindent is called \emph{evolution cocycle} over the evolution
semiflow $\psi$.
\end{definition}

\begin{definition}\rm\label{d_lses}
A function $\xi:\mathcal{T}\times \mathcal{Y}\rightarrow
\mathcal{Y}$ defined by
\begin {equation}
\xi(t,s,x,v)=(\psi(t,s,x),\Psi(t,s,x)v), \ \forall (t,s,x,v)\in
\mathcal{T}\times \mathcal{Y}
\end{equation}
where $\Psi$ is an evolution cocycle over the evolution semiflow
$\psi$, is called \emph{skew-evolution semiflow} on $\mathcal{Y}$.
\end{definition}

\begin{example}\rm\label{ex_cev}
Let $f:\mathbb{R}\rightarrow\mathbb{R}_{+}$ be a function which is
nondecreasing on the interval $(-\infty,0)$ and decreasing on the
interval $(0,\infty)$ with the property that there exists
\[
\underset{t\rightarrow\pm\infty}{\lim }f(t)\in(0,\infty).
\]
Let us consider the set
$\mathcal{C}=\mathcal{C}(\mathbb{R},\mathbb{R})$ of all continuous
functions given by $x:\mathbb{R}\rightarrow \mathbb{R}$, endowed
with the uniform convergence topology on compact subsets of
$\mathbb{R}$. $\mathcal{C}$ is metrizable by respect to the metric
\begin{equation*}
d(x,y)=\sum_{n=1}^{\infty}\frac{1}{2^{n}}\frac{d_{n}(x,y)}{1+d_{n}(x,y)}
\end{equation*}
where
\begin{equation*}
d_{n}(x,y)=\underset{t\in [-n,n]}\sup{|x(t)-y(t)|}.
\end{equation*}

Let $\mathcal{X}$ be the closure in $\mathcal{C}$ of the set ${\{f_{t}, \ t\in \mathbb{R}\}}$, where
\[
f_{t}(\tau)=f(t+\tau), \ \forall \tau\in \mathbb{R}, f\in
\mathcal{C}.
\]
Then $(\mathcal{X},d)$ is a metric space and the
mapping
\begin{equation*}
\psi:\mathcal{T}\times \mathcal{X}\rightarrow \mathcal{X}, \
\psi(t,s,x)=x_{t-s}
\end{equation*}
is an evolution semiflow on $\mathcal{X}$.

The mapping
\begin{equation*}
\Psi:\mathcal{T}\times \mathcal{X}\rightarrow
\mathfrak{B}(\mathcal{V}),
\end{equation*}
given by
\begin{equation*}
\Psi(t,s,x)(v)=e^{\int_{s}^{t}x(\tau-s)d\tau}v
\end{equation*}
is an evolution cocycle.

Then $\xi=(\psi,\Psi)$ is a skew-evolution semiflow on
$\mathcal{Y}$.
\end{example}

\begin{definition}\rm
A mapping $P:\mathcal{Y}\rightarrow \mathcal{Y}$ is said to be a
\emph{projector} on $\mathcal{Y}$ if $P$ is continuous and has the
form
\begin{equation}
P(x,v)=(x,P(x)v), \ (x,v)\in \mathcal{Y}
\end{equation}
where $P(x)$ is a projection on $\mathcal{Y}_{x}=\{x\}\times
\mathcal{V}$, $x\in \mathcal{X}$.
\end{definition}

\begin{remark}\rm
The function $P(x):\mathcal{Y}_{x}\rightarrow \mathcal{Y}_{x}$ is
a bounded mapping with the property $P(x)P(x)=P^{2}(x)=P(x)$ for
all $x\in \mathcal{X}.$
\end{remark}

\begin{definition}\rm
The mapping $Q:\mathcal{Y}\rightarrow \mathcal{Y}$ given by
\begin{equation}
Q(x,v)=(x,v-P(x)v)
\end{equation}
where $P$ is a projector on $\mathcal{Y}$, is called
\emph{complementary projector} to $P$ on $\mathcal{Y}$.
\end{definition}

\begin{definition}\rm\label{projector}
A projector $P$ on $\mathcal{Y}$ is said to be invariant on
$\mathcal{Y}$ relative to a skew-evolution semiflow $\xi=(
\psi,\Psi)$ if one has
\begin{equation}
P(\psi(t,s,x))\Psi(t,s,x)v=\Psi(t,s,x)P(x)v,
\end{equation}
\noindent for all $(t,s)\in \mathcal{T}$ and for all $(x,v)\in
\mathcal{Y}$, where $\Psi$ denotes an evolution cocycle over the
semiflow $\psi$.
\end{definition}

\begin{remark}\rm
If the projector $P$ is invariant, then the complementary
projector $Q$ is invariant as well.
\end{remark}

Let us denote $\Psi_{k}(t,t_{0},x)=\Psi(t,t_{0},x)P_{k}(x)$, for
all $(t,t_{0})\in \mathcal{T}$, all $x\in \mathcal{X}$ and $k\in
\{0,1,2\}.$

\begin{definition}\rm
Three projector families $\{P_{k}(y)\}_{y\in \mathcal{Y}}$, $k\in
\{0,1,2\}$, are said to be \emph{compatible} with a skew-evolution
semiflow $\xi$ if

$(c_{1})$ each of the projectors $P_{k}$, $k\in \{0,1,2\}$ is
invariant on $\mathcal{Y}$

$(c_{2})$ for each $x\in \mathcal{X}$ and for all  $i,j \in
\{0,1,2\}, \ i\neq j$
\begin{equation}
P_{0}(x)+P_{1}(x)+P_{2}(x)=I \ \textrm{and} \ P_{i}(x)P_{j}(x)=0.
\end{equation}
\end{definition}

\begin{definition}\rm\label{uet}
A skew-evolution
semiflow $\xi=(\psi,\Psi)$ has \emph{uniform exponential trichotomy} on $\mathcal{Y}$ if there exist $%
N_{0},N_{1},N_{2}>1$, $\nu _{0},\nu _{1},\nu _{2}>0$ and three
projector families $\{P_{k}(y)\}_{y\in \mathcal{Y}}$, $k\in
\{0,1,2\}$, compatible with $\xi$ such that

$(t_{0})$
\begin{equation}
\left\Vert P_{0}(x)v\right\Vert \leq N_{0}e^{\nu
_{0}(t-t_{0})}\left\Vert \Psi_{0}(t,t_{0},x)v\right\Vert \leq
N_{0}^{2}e^{2\nu _{0}(t-t_{0})}\left\Vert P_{0}(x)v\right\Vert
\end{equation}

$(t_{1})$
\begin{equation}
e^{\nu _{1}(t-t_{0})}\left\Vert \Psi_{1}(t,t_{0},x)v\right\Vert
\leq N_{1}\left\Vert P_{1}(x)v\right\Vert
\end{equation}

$(t_{2})$
\begin{equation}
e^{\nu _{2}(t-t_{0})}\left\Vert P_{2}(x)v\right\Vert \leq
N_{2}\left\Vert \Psi_{2}(t,t_{0},x)v\right\Vert
\end{equation}

\noindent for all $(t,t_{0})\in \mathcal{T}$ and all $(x,v)\in
\mathcal{Y}.$
\end{definition}

\begin{remark}\rm
$(i)$ For $P_{0}=0$ in Definition \ref{uet} we obtain the property
of uniform exponential dichotomy.

$(ii)$ If $P_{0}=P_{2}=0$ the property of uniform exponential
stability is obtained. It follows that a uniformly exponentially
stable skew-evolution semiflow is uniformly exponentially
dichotomic and, further, uniformly exponentially trichotomic.

$(iii)$ Also it is easy to observe that the property of uniform
exponential instability implies the uniform exponential dichotomy
and, further, the uniform exponential trichotomy.
\end{remark}

\begin{example}\rm\label{ex_uet}
Let $f:\mathbb{R}_{+}\rightarrow (0, \infty)$ be a decreasing
function with the property that there exists
$\underset{t\rightarrow \infty}\lim{f(t)}=l>0$. Let $\mu>f(0)>0$.

Let $(\mathcal{X},d)$ be the metric space given as in Example
\ref{ex_cev}, the closure in $\mathcal{C}(\mathbb{R},\mathbb{R})$
of the set ${\{f_{t}, \ t\in \mathbb{R}\}}$, where
\[
f_{t}(\tau)=f(t+\tau), \ \forall \tau\in \mathbb{R}, f\in
\mathcal{C}.
\]
The mapping
\begin{equation*}
\psi:\mathcal{T}\times \mathcal{X}\rightarrow \mathcal{X}, \
\psi(t,s,x)=x_{t-s},
\end{equation*}
is an evolution semiflow on $\mathcal{X}$.

We consider $\mathcal{V}=\mathbb{R}^{3}$ with the norm
\begin{equation*}
\left\Vert v\right\Vert=|v_{1}|+|v_{2}|+|v_{3}|, \
v=(v_{1},v_{2},v_{3})\in \mathcal{\mathbb{R}}^{3}.
\end{equation*}

The mapping $\Psi:\mathcal{T}\times \mathcal{X} \rightarrow
\mathcal{B}(\mathcal{\mathbb{R}}^{3})$ defined by
\[
\Psi(t,t_{0},x)v=
\]
\[
=(e^{-\mu(t-t_{0})+\int_{t_{0}}^{t}x(\tau-t_{0})d\tau}v_{1},\
e^{\int_{t_{0}}^{t}x(\tau-t_{0})d\tau}v_{2},\
e^{-(t-t_{0})x(0)+\int_{t_{0}}^{t}x(\tau-t_{0})d\tau}v_{3}),
\]
for all $(t,t_{0})\in \mathcal{T}$ and all $(x,v)\in \mathcal{Y}$,
is an evolution cocycle.

We consider the projections given by
\[
P_{1}(x)(v)=(v_{1},0,0), \ P_{2}(x)(v)=(0,v_{2},0), \
P_{0}(x)(v)=(0,0,v_{3})
\]
for all $x\in \mathcal{X}$ and all $v=(v_{1},v_{2},v_{3})\in
\mathcal{V}$.

We obtain the relations
\[
\left\Vert \Psi(t,t_{0},x)P_{1}(x)v)\right\Vert \leq
e^{[-\mu+x(0)](t-s)}\left\Vert
\Psi(s,t_{0},x)P_{1}(x)v)\right\Vert
\]
\[
\left\Vert \Psi(t,t_{0},x)P_{2}(x)v)\right\Vert \geq
e^{l(t-s)}\left\Vert \Psi(s,t_{0},x)P_{2}(x)v)\right\Vert
\]
\[
\left\Vert \Psi(t,t_{0},x)P_{0}(x)v)\right\Vert \leq
e^{x(0)(t-s)}\left\Vert \Psi(s,t_{0},x)P_{0}(x)v)\right\Vert
\]
\[
\left\Vert \Psi(t,t_{0},x)P_{0}(x)v)\right\Vert \geq
e^{-x(0)(t-s)}\left\Vert \Psi(s,t_{0},x)P_{0}(x)v)\right\Vert
\]
for all $(t,s),(s,t_{0})\in \mathcal{T}$ and all $(x,v)\in
\mathcal{Y}$.

It follows that the skew-evolution semiflow given by $\xi=(\psi,
\Psi)$ is uniformly exponentially trichotomic with characteristics
\[
N_{0}=N_{1}=N_{2}=1 \ \textrm{and} \ \nu_{0} =x(0), \
\nu_{1}=-\mu+x(0), \ \nu_{2}=l.
\]
\end{example}

\vspace{3mm}

In what follows by projections we will consider the mappings
$\widetilde{P}:\mathcal{Y}\rightarrow \mathcal{Y}$ with the
property $\widetilde{P}^{2}(t)=\widetilde{P}(t)$, $\forall t\geq
0$.

\begin{definition}\rm\label{projection}
A family of projections $\{\widetilde{P}(t)\}_{t\geq 0}$ is said
to be \emph{invariant} on $\mathcal{Y}$ relative to a
skew-evolution semiflow $\xi=( \psi,\Psi)$ if one has
\begin{equation}\label{rel_comp_pct}
\widetilde{P}(t+s)\Psi(t,s,\psi(t,s,x))v=\Psi(t,s,\psi(t,s,x))\widetilde{P}(s)v,
\end{equation}
\noindent for all $(t,s)\in \mathcal{T}$ and for all $(x,v)\in
\mathcal{Y}$, where $\Psi$ denotes an evolution cocycle over the
semiflow $\psi$.
\end{definition}

\begin{definition}\rm
Three projection families $\{\widetilde{P}_{k}(t)\}_{t\geq 0}$,
$k\in \{0,1,2\}$, are said to be \emph{compatible} with a
skew-evolution semiflow $\xi$ at point $x\in \mathcal{X}$ if

$(cp_{1})$ each projection $\widetilde{P}_{k}$, $k\in \{0,1,2\}$
is invariant on $\mathcal{Y}$

$(cp_{2})$ for each $t\in \mathbb{R}_{+}$ one has
\begin{equation*}
\widetilde{P}_{0}(t)+\widetilde{P}_{1}(t)+\widetilde{P}_{2}(t)=I \
\textrm{and} \ \widetilde{P}_{i}(t)\widetilde{P}_{j}(t)=0, \
\forall i,j \in \{0,1,2\}, \ i\neq j
\end{equation*}

$(cp_{3})$ for all $t\geq 0$, all $v\in \mathcal{V}$ and all
$i,j\in \{0,1,2\},\ i\neq j$
\[
\left\Vert
\widetilde{P}_{i}(t)v+\widetilde{P}_{j}(t)v\right\Vert^{2}=\left\Vert
\widetilde{P}_{i}(t)v\right\Vert^{2}+\left\Vert
\widetilde{P}_{j}(t)v\right\Vert^{2}.
\]
\end{definition}

\begin{definition}\rm\label{puet}
A skew-evolution semiflow $\xi=(\psi,\Psi)$ has \emph{pointwise
uniform exponential trichotomy}
on $\mathcal{Y}$ at point $x_{0}\in \mathcal{X}$ if there exist some constants $%
N_{0}^{x_{0}},N_{1}^{x_{0}},N_{2}^{x_{0}}>1$, $\nu
_{0}^{x_{0}},\nu _{1}^{x_{0}},\nu _{2}^{x_{0}}>0$ and three
projection families $\{\widetilde{P}_{k}(t)\}_{t\geq 0}$, $k\in
\{0,1,2\}$, compatible with $\xi$ such that

$(pt_{0})$
\[
\left\Vert
\Psi(s,t_{0},\psi(s,t_{0},x_{0})\widetilde{P}_{0}(t_{0})v\right\Vert
\leq
\]
\begin{equation}\label{pt01}
\leq N_{0}^{x_{0}}e^{\nu _{0}^{x_{0}}(t-s)}\left\Vert
\Psi(t,t_{0},\psi(t,t_{0},x_{0}))\widetilde{P}_{0}(t_{0})v\right\Vert
\end{equation}
\[
\left\Vert\Psi(t,t_{0},\psi(t,t_{0},x_{0}))\widetilde{P}_{0}(t_{0})v\right\Vert
\leq
\]
\begin{equation}\label{pt02}
\leq N_{0}^{x_{0}}e^{\nu _{0}^{x_{0}}(t-s)}\left\Vert
\Psi(s,t_{0},\psi(s,t_{0},x_{0}))\widetilde{P}_{0}(t_{0})v\right\Vert
\end{equation}

$(pt_{1})$
\[
e^{\nu _{1}^{x_{0}}(t-s)}\left\Vert
\Psi(t,t_{0},\psi(t,t_{0},x_{0}))\widetilde{P}_{1}(t_{0})v\right\Vert
\leq
\]
\begin{equation}\label{pt1}
\leq N_{1}^{x_{0}}\left\Vert
\Psi(s,t_{0},\psi(s,t_{0},x_{0}))\widetilde{P}_{1}(t_{0})v\right\Vert
\end{equation}

$(pt_{2})$
\[
e^{\nu _{2}^{x_{0}}(t-s)}\left\Vert
\Psi(s,t_{0},\psi(s,t_{0},x_{0}))\widetilde{P}_{2}(t_{0})v\right\Vert
\leq
\]
\begin{equation}\label{pt2}
\leq N_{2}^{x_{0}}\left\Vert
\Psi(t,t_{0},\psi(t,t_{0},x_{0}))\widetilde{P}_{2}(t_{0})v\right\Vert
\end{equation}

\noindent for all $(t,s),(s,t_{0})\in \mathcal{T}$ and all $v\in
\mathcal{V}.$
\end{definition}

\begin{remark}\rm
According to \cite{LaMoRa_JDE} the concept introduced in
Definition \ref{puet} can also be called exponential trichotomy of
the skew-evolution semiflow $\xi$ over the orbit through $x_{0}\in
\mathcal{X}$.
\end{remark}

\begin{remark}\rm
Referring to the case of pointwise dichotomy, in \cite{LaMoRa_JDE}
and \cite{SaSe_JDE} is shown that for an invertible-valued cocycle
involved in the definition of a skew-product flow, the later
property is equivalent with the classical definition of global
dichotomy.

We remark that the property of uniform exponential trichotomy for
a skew-evolution semiflow implies the uniform exponential
trichotomy at each point $x\in \mathcal{X}$, by considering
$\widetilde{P}_{k}(t)=P_{k}(\psi(t,x)), \ \forall t \geq 0$, $k\in
\{0,1,2\}$.

Instead, in general, the converse affirmation is not always true,
as shown in the following
\end{remark}

\begin{example}\rm\label{contrex_uet}
Let
$\mathcal{C}(\mathbb{R}_{+},\mathbb{R})=\{f:\mathbb{R}_{+}\rightarrow
\mathbb{R} \ | \ f \ \textrm{continuous}\}$, metrizable with
respect to the metric considered in Example \ref{ex_cev}.

For every $n\in \mathbb{N}^{*}$ we consider a decreasing function
\[
x_{n}:\mathbb{R}_{+}\rightarrow
\left(\frac{1}{2n+1},\frac{1}{2n}\right),\ \underset{t\rightarrow
\infty}\lim{x_{n}(t)}=\frac{1}{2n+1}.
\]
If $x_{n}^{s}(t)=x_{n}(t+s), \ \forall t,s \geq 0$ and
$\mathcal{X}=\overline{\{x_{n}^{s},n\in \mathbb{N}^{*},s\in
\mathbb{R}_{+}\}}$ then the mapping
\begin{equation*}
\psi:\mathcal{T}\times \mathcal{X}\rightarrow \mathcal{X}, \
\psi(t,s,x)(\tau)=x(t-s+\tau)
\end{equation*}
is an evolution semiflow on $\mathcal{X}$.

If we consider $\mathcal{V}=\mathbb{R}^{3}$ with the norm
$\left\Vert
(v_{1},v_{2},v_{3})\right\Vert=|v_{1}|+|v_{2}|+|v_{3}|$, then
$\xi=(\psi,\Psi)$, where $\Psi:\mathcal{T}\times
\mathcal{X}\rightarrow \mathfrak{B}(\mathcal{V})$ is given by
\begin{equation*}
\Psi(t,t_{0},x)(v_{1},v_{2},v_{3})=\left(e^{-\int_{t_{0}}^{t}x(s)ds}v_{1},e^{\int_{t_{0}}^{t}x(s)ds}v_{2},e^{\int_{t_{0}}^{t}x(s)ds}v_{3}\right),
\end{equation*}
is a skew-evolution semiflow on $\mathcal{Y}$. We consider the
projections $P_{k}$, $k\in \{0,1,2\}$ given as in Example
\ref{ex_uet}. We obtain that $\xi$ is not globally uniformly
exponentially trichotomic, but is uniformly exponentially
trichotomic at every point $x\in \mathcal{X}$ relative to the
projection families $\{\widetilde{P}_{k}(x)\}_{x\in \mathcal{X}}$,
given by $\widetilde{P}_{k}(x)=P_{k}$, $x\in \mathcal{X}$, $k\in
\{0,1,2\}$.
\end{example}

\section{Characterizations of the pointwise trichotomy}

Let us consider the set of functions
\[
\mathcal{F}=\{f:[0,\infty )\rightarrow (0,\infty )\ | \ f \
\textrm{decreasing function}, \ \underset{t\rightarrow \infty
}{\lim }f(t)=0\}
\]
In order to emphasize the natural extension of the trichotomy
relative to the property of dichotomy, we will present a
characterization by means of two projection families compatible
with a skew-evolution semiflow, introduced by the following

\begin{definition}\label{comp2}\rm
Two projection families $\{\widetilde{Q}_{k}(t)\}_{t\geq 0}$ where
$\widetilde{Q}_{k}:\mathbb{R}_{+}\rightarrow
\mathfrak{B}(\mathcal{V})$, $k\in \{1,2\}$, are said to be
\textit{compatible} with a skew-evolution semiflow $\xi$ at point
$x\in \mathcal{X}$ if following relations hold

$(cq_{1})$ $%
\widetilde{Q}_{1}(t)\widetilde{Q}_{2}(t)=\widetilde{Q}_{2}(t)\widetilde{Q}_{1}(t)=0$

$(cq_{2})$ $\left\Vert \left[ \widetilde{Q}_{1}(t)+\widetilde{Q}_{2}(t)%
\right] v\right\Vert ^{2}=\left\Vert
\widetilde{Q}_{1}(t)v\right\Vert ^{2}+\left\Vert
\widetilde{Q}_{2}(t)v\right\Vert ^{2}$

$(cq_{3})$ $\left\Vert \left[ I-\widetilde{Q}_{1}(t)\right]
v\right\Vert ^{2}=\left\Vert \left[
I-\widetilde{Q}_{1}(t)-\widetilde{Q}_{2}(t)\right] v\right\Vert
^{2}+\left\Vert \widetilde{Q}_{2}(t)v\right\Vert ^{2}$

$(cq_{4})$ $\left\Vert \left[ I-\widetilde{Q}_{2}(t)\right]
v\right\Vert ^{2}=\left\Vert \left[
I-\widetilde{Q}_{1}(t)-\widetilde{Q}_{2}(t)\right] v\right\Vert
^{2}+\left\Vert \widetilde{Q}_{1}(t)v\right\Vert ^{2}$

$(cq_{5})$ $%
\Psi(t,t_{0},\psi(t,t_{0},x))\widetilde{Q}_{k}(t_{0})v=\widetilde{Q}_{k}(t+t_{0})\Psi(t,t_{0},\psi(t,t_{0},x))v$

\noindent for all $(t,t_{0})\in \mathcal{T}$, all $v\in
\mathcal{V}$ and $k\in \{1,2\}$.
\end{definition}

\begin{theorem}
A skew-evolution semiflow $\xi=(\psi,\Psi)$ is pointwise uniformly
exponentially trichotomic at point $x_{0}\in \mathcal{X}$ if and
only if there exist two functions $\varphi_{1},\varphi_{2} \in
\mathcal{F}$ and two projection families
$\{\widetilde{Q}_{k}(t)\}_{t\geq 0}$, $k\in \{1,2\}$, compatible
with $\xi$ such that following inequalities hold

$(puet_{0})$ $\varphi_{1}(t-s)\left\Vert
[I-\widetilde{Q}_{1}(s)]v\right\Vert \leq \left\Vert
[I-\widetilde{Q}_{1}(t+s)]\Psi(t,s,\psi(t,s,x_{0})) v\right\Vert$

$(puet_{0}')$ $\varphi_{1}(t-s)\left\Vert
[I-\widetilde{Q}_{2}(t+s)] \Psi(t,s,\psi(t,s,x_{0}))v\right\Vert
\leq \left\Vert [I-\widetilde{Q}_{2}(s)]v\right\Vert$

$(puet_{1})$ $\left\Vert
\Psi(t,s,\psi(t,s,x_{0}))\widetilde{Q}_{1}(s)v\right\Vert \leq
\varphi_{2}(t-s)\left\Vert \widetilde{Q}_{1}(s)v\right\Vert$

$(puet_{2})$ $\left\Vert
\widetilde{Q}_{2}(s)v\right\Vert\leq\varphi_{2}(t-s)\left\Vert
\Psi(t,s,\psi(t,s,x_{0}))\widetilde{Q}_{2}(s)v\right\Vert$

\noindent for all $(t,s)\in \mathcal{T}$ and all $v\in
\mathcal{V}$.
\end{theorem}

\begin{proof}
\emph{Necessity.} The existence of functions
$\varphi_{1},\varphi_{2}\in \mathcal{F}$ and projections
$\widetilde{Q}_{1}$ and $\widetilde{Q}_{2}$ is assured by
Definition \ref{puet} if we consider
\[
\varphi_{1}(t)=N_{0}^{x_{0}}e^{-\nu_{0}^{x_{0}}t} \ \textrm{and} \
\varphi_{2}(t)=\max\{N_{1}^{x_{0}},N_{2}^{x_{0}}\}e^{-\min\{\nu_{1}^{x_{0}},\nu_{2}^{x_{0}}\}t},
\ t\geq 0
\]
respectively
\[
\widetilde{Q_{1}}(t)=\widetilde{P_{1}}(t), \
\widetilde{Q_{2}}(t)=\widetilde{P_{2}}(t), \ \ t\geq 0
\]
and which verify the compatibility (see \cite{MeSt_OT}) and
$(puet_{0})$--$(puet_{2})$.

\emph{Sufficiency.} We define
\[
\widetilde{P}_{0}(t)=I-\widetilde{Q}_{1}(t)-\widetilde{Q}_{2}(t),
\ \widetilde{P}_{1}(t)=\widetilde{Q}_{1}(t),\
\widetilde{P}_{2}(t)=\widetilde{Q}_{2}(t), \ t\geq 0.
\]
The compatibility of the projections families
$\{\widetilde{P}_{k}(t)\}_{t\geq 0}$, $k\in \{0,1,2\}$ with $\xi$
is easy to verify (see \cite{MeSt_OT}). By the properties of
$\widetilde{Q}_{1}$ and $\widetilde{Q}_{2}$, it is true that
\[
\widetilde{P}_{0}(t)=[I-\widetilde{Q}_{1}(t)][I-\widetilde{Q}_{2}(t)],
\ \forall t\geq 0.
\]

Let $x_{0}\in \mathcal{X}$ and $\varphi_{1}\in \mathcal{F}$. To
prove the existence of $N_{0}^{x_{0}}$ and $\nu_{0}^{x_{0}}$
required by relation (\ref{pt01}) we consider following
inequalities which hold for all $(t,t_{0})\in \mathcal{T}$ and all
$v\in \mathcal{V}$
\[
\left\Vert
\widetilde{P}_{0}(t+t_{0})\Psi(t,t_{0},\psi(t,t_{0},x_{0}))v\right\Vert\geq
\]
\[
\geq \varphi_{1}(1)\left\Vert
\Psi(t-1,t_{0},\psi(t-1,t_{0},x_{0}))[I-\widetilde{Q}_{1}(t_{0})][I-\widetilde{Q}_{2}(t_{0})]v\right\Vert\geq
...\geq
\]
\[
\geq \varphi_{1}^{[t-t_{0}]}(1)\left\Vert
\Psi(t-[t-t_{0}],t_{0},\psi(t_{0},t_{0},x_{0}))\widetilde{P}_{0}(t_{0})v\right\Vert\geq
\]
\[
\geq\varphi_{1}^{[t-t_{0}]}(\delta)\varphi_{1}(1)\left\Vert
\Psi(t_{0},t_{0},\psi(t_{0},t_{0},x_{0}))\widetilde{P}_{0}(t_{0})v\right\Vert,
\]
where by $[t-t_{0}]$ we have denoted the integer part of the
considered difference and the existence of $\delta>1$ is assured
by the definition of the function set $\mathcal{F}$, such that
$\varphi_{1}(\delta)<1$.

If we denote $N_{0}^{x_{0}}=\varphi_{1}(1)$ and
$\nu_{0}^{x_{0}}=-\ln \varphi_{1}(\delta)$, the relation
(\ref{pt01}) is proved. Analogues techniques are used to prove
inequality (\ref{pt02}), as well as (\ref{pt1}) and (\ref{pt2}).
\end{proof}

\vspace{3mm}

Some integral characterizations for the pointwise trichotomy are
given by the next result, involving four compatible projector
families, given by the next conditions.

\begin{definition}\label{comp4}\rm
Four projection families $\{\widetilde{R}_{k}(t)\}_{t\geq 0}$
where $\widetilde{R}_{k}:\mathbb{R}_{+}\rightarrow
\mathfrak{B}(\mathcal{V})$, $k\in \{1,2,3,4\}$, are said to be
\textit{compatible} with a skew-evolution semiflow
$\xi=(\psi,\Psi)$ at point $x\in \mathcal{X}$ if

$(cr_{1})$ $%
\widetilde{R}_{1}(t)+\widetilde{R}_{3}(t)=\widetilde{R}_{2}(t)+\widetilde{R}_{4}(t)=I$

$(cr_{2})$ $%
\widetilde{R}_{1}(t)\widetilde{R}_{2}(t)=\widetilde{R}_{2}(t)\widetilde{R}_{1}(t)=0$ and $%
\widetilde{R}_{3}(t)\widetilde{R}_{4}(t)=\widetilde{R}_{4}(t)\widetilde{R}_{3}(t)$

$(cr_{3})$ $\left\Vert \left[
\widetilde{R}_{1}(t)+\widetilde{R}_{2}(t)\right] v\right\Vert
^{2}=\left\Vert \widetilde{R}_{1}(t)v\right\Vert ^{2}+\left\Vert
\widetilde{R}_{2}(t)v\right\Vert ^{2}$

$(cr_{4})$ $\left\Vert \left[
\widetilde{R}_{1}(t)+\widetilde{R}_{3}(t)\widetilde{R}_{4}(t)\right]
v\right\Vert ^{2}=\left\Vert \widetilde{R}_{1}(t)v\right\Vert
^{2}+\left\Vert
\widetilde{R}_{3}(t)\widetilde{R}_{4}(t)v\right\Vert ^{2}$

$(cr_{5})$ $\left\Vert \left[
\widetilde{R}_{2}(t)+\widetilde{R}_{3}(t)\widetilde{R}_{4}(t)\right]
v\right\Vert ^{2}=\left\Vert \widetilde{R}_{2}(t)v\right\Vert
^{2}+\left\Vert
\widetilde{R}_{3}(t)\widetilde{R}_{4}(t)v\right\Vert ^{2}$

$(cr_{6})$ $%
\Psi(t,t_{0},\psi(t,t_{0},x))\widetilde{R}_{k}(t_{0})v=\widetilde{R}_{k}(t+t_{0})\Psi(t,t_{0},\psi(t,t_{0},x))v$

\noindent for all $t,t_{0}\geq 0, \ t\geq t_{0}$, for all
$v\in\mathcal{V}$ and $k\in \{1,2,3,4\}$
\end{definition}

\begin{theorem}
Let $\xi=(\psi,\Psi)$ be a skew-evolution semiflow such that for
each $(x,v)\in \mathcal{Y}$ the mapping $s\rightarrow
\Phi(s,t_{0},x)v $ is continuous on $[t_{0},\infty )$ and the
mapping $s\rightarrow \Phi(t,s,x)v$ is continuous on $[0,t]$ and
let us consider four projection families
$\{\widetilde{R}_{k}(t)\}_{t\geq 0}$, $k\in \{1,2,3,4\}$,
compatible with $\xi=(\psi,\Psi)$, with the property that there
exist $s_{0}>t_{0}\geq 0$ and $c\in (0,1)$ such that
\begin{equation}\label{st}
\left\Vert
\Psi(s_{0},t_{0},\psi(s_{0},t_{0},x))\widetilde{R}_{1}(t_{0})v\right\Vert
\leq c\left\Vert \widetilde{R}_{1}(t_{0})v\right\Vert
\end{equation}%
\begin{equation}\label{in}
\left\Vert \widetilde{R}_{2}(t_{0})v\right\Vert \leq c\left\Vert
\Psi(s_{0},t_{0},\psi(s_{0},t_{0},x))\widetilde{R}_{2}(t_{0})v\right\Vert
\end{equation}%
\noindent for all $(x,v)\in \mathcal{Y}$.

Let $N, \widetilde{N}>1$ and $\omega,\widetilde{\omega}>0$ such
that
\begin{equation}\label{eg}
\left\Vert \Psi(t,t_{0},x)\widetilde{R}_{i}(t_{0})v\right\Vert
\leq Ne^{\omega(t-s)}\left\Vert
\Psi(s,t_{0},x)\widetilde{R}_{i}(t_{0})v\right\Vert, \ i\in\{1,3\}
\end{equation}%
\begin{equation}\label{ed}
\left\Vert \Psi(s,t_{0},x)\widetilde{R}_{j}(t_{0})v\right\Vert
\leq \widetilde{N}e^{\widetilde{\omega}(t-s)}\left\Vert
\Psi(t,t_{0},x)\widetilde{R}_{j}(t_{0})v\right\Vert, \ j\in
\{2,4\}
\end{equation}%
\noindent for all $(t,s),(s,t_{0})\in\mathcal{T}$ and all
$(x,v)\in \mathcal{Y}$.

Then the skew-product semiflow is pointwise uniformly
exponentially trichotomic at point $x_{0}\in \mathcal{X}$ if and
only if

$(upet_{0})$ there exist $M>0$ and $%
\alpha >0$ such that
\begin{equation*}
\int_{t_{0}}^{t}e^{-\alpha (\tau -t_{0})}\left\Vert
\Psi(\tau,t_{0},\psi(\tau,t_{0},x_{0})
)\widetilde{R}_{3}(t_{0})v\right\Vert d\tau \leq M\left\Vert
\widetilde{R}_{3}(t_{0})v\right\Vert
\end{equation*}%
for all $(t,t_{0})\in \mathcal{T}$ and all $v\in \mathcal{V}$

$(upet_{0}')$ there exist $M>0$ and $%
\alpha >0$ such that
\begin{equation*}
\int_{s}^{t}e^{\alpha (t-\tau )}\left\Vert
\Psi(\tau,t_{0},\psi(\tau,t_{0},x_{0})
)\widetilde{R}_{4}(t_{0})v\right\Vert d\tau \leq M\left\Vert
\Psi(t,t_{0},\psi(t,t_{0},x_{0}))\widetilde{R}_{4}(t_{0})v\right\Vert
\end{equation*}%
for all $(t,s),(s,t_{0})\in \mathcal{T}$ and all $v\in V$

$(upet_{1})$ there exists $M_{1}>0$ such that
\begin{equation*}
\int_{t_{0}}^{\infty }\left\Vert
\Psi(\tau,t_{0},\psi(\tau,t_{0},x_{0})
)\widetilde{R}_{1}(t_{0})v\right\Vert d\tau \leq M_{1}\left\Vert
\widetilde{R}_{1}(t_{0})v\right\Vert
\end{equation*}%
for all $t_{0}\in \mathbb{R}_{+}$ and all $v\in V$

$(upet_{2})$ there exists $M_{2}>0$ such that
\begin{equation*}
\int_{0}^{t}\left\Vert
\Psi(\tau,t_{0},\psi(\tau,t_{0},x_{0}))\widetilde{R}_{2}(t_{0})v\right\Vert
d\tau \leq M_{2}\left\Vert
\Psi(t,t_{0},\psi(t,t_{0},x_{0}))\widetilde{R}_{2}(t_{0})v\right\Vert
\end{equation*}%
for all $(t,t_{0})\in \mathcal{T}$ and all $v\in V$.
\end{theorem}

\begin{proof} \emph{Necessity}. Considering the projection families
$\{\widetilde{P}_{k}(t)\}_{t\geq 0}$, $k\in \{0,1,2\}$, given as
in Definition \ref{puet} we define
\[
\widetilde{R}_{1}(t)=\widetilde{P}_{1}(t), \
\widetilde{R}_{2}(t)=\widetilde{P}_{2}(t), \
\widetilde{R}_{3}(t)=I-\widetilde{P}_{1}(t), \
\widetilde{R}_{4}(t)=I-\widetilde{P}_{2}(t), \ t\geq 0.
\]
We also have
$\widetilde{R}_{3}(t)\widetilde{R}_{4}(t)=\widetilde{R}_{4}(t)\widetilde{R}_{3}(t)=\widetilde{P}_{0}(t),
\ t\geq 0$. The compatibility is easy to verify, as in
\cite{MeSt_OT}. We denote by $N^{x_{0}}=\max
\{N_{0}^{x_{0}},N_{1}^{x_{0}},N_{2}^{x_{0}}\}$ and
$\nu_{1}^{x_{0}},\nu_{2}^{x_{0}},\nu_{0}^{x_{0}}$ the
characteristics of the pointwise uniformly exponentially
trichotomic skew-evolution semiflow $\xi$, given by Definition
\ref{puet}.

To prove $(upet_{1})$, we consider $M_{1}=N^{x_{0}}(\nu
_{1}^{x_{0}})^{-1}$. For $(upet_{2})$ let us denote
$M_{2}=N^{x_{0}}(\nu _{2}^{x_{0}})^{-1}$. We obtain $(upet_{0})$
and $(upet_{0}')$ if we consider $\alpha=2\nu_{0}^{x_{0}}$ and
$M=N^{x_{0}}(\nu_{0}^{x_{0}})^{-1}$.

\emph{Sufficiency}. We define
\[
\widetilde{P}_{1}(t)=\widetilde{R}_{1}(t), \
\widetilde{P}_{2}(t)=\widetilde{R}_{2}(t), \
\widetilde{P}_{0}(t)=\widetilde{R}_{3}(t)\widetilde{R}_{4}(t), \
t\geq 0.
\]
The compatibility is obtained from the properties of the
projection families $\{\widetilde{R}_{k}(t)\}_{t\geq 0}$, $k\in
\{1,2,3,4\}$.

Relation $(upet_{0})$ defines the integral stability for the
skew-evolution flow given by $\xi_{\alpha}=(\psi,\Psi_{\alpha})$,
where $\Psi_{\alpha}(t,s,x)=e^{-\alpha (t-s)}\Psi(t,s,x)$,
$(t,s)\in\mathcal{T}, \ x\in \mathcal{X}$. As relation (\ref{eg})
hold, then according to a result that characterize the stability
obtained in \cite{MeStBu_UVT} it follows that $\xi_{\alpha}$ is
uniformly exponentially stable, which assures the existence of
some constants $K>0$ and $\beta>0$ such that we obtain
successively for all $(t,s),(s,t_{0})\in\mathcal{T}$ and all
$(x,v)\in \mathcal{Y}$
\[
\left\Vert
\Psi(t,t_{0},\psi(t,t_{0},x))\widetilde{R}_{3}(t_{0})v)\right\Vert
=e^{\alpha(t-t_{0})}\left\Vert
\Psi_{\alpha}(t,t_{0},\psi(t,t_{0},x))\widetilde{R}_{3}(t_{0})v)\right\Vert\leq
\]
\[
\leq Ke^{-\beta(t-s)}e^{\alpha(t-t_{0})}\left\Vert
\Psi_{\alpha}(s,t_{0},\psi(s,t_{0},x))\widetilde{R}_{3}(t_{0})v)\right\Vert=
\]
\[
=Ke^{(\alpha-\beta)(t-s)}\left\Vert
\Psi(s,t_{0},\psi(s,t_{0},x))\widetilde{R}_{3}(t_{0})v)\right\Vert
\]
and, if we define
\begin{equation*}
\nu_{0}^{x_{0}} =\left\{
\begin{array}{cc}
\alpha -\beta , & \text{if }\alpha >\beta \\
1, & \text{if }\alpha \leq \beta%
\end{array}%
\right.
\end{equation*}

\noindent then (\ref{pt01}) is proved. By an analogue proof is
obtained (\ref{pt02}) from $(upet_{0}')$.

As in \cite{MeSaSa_MB}, from the property that there exist
$t_{0}>0$ and $c\in (0,1)$ such that relation (\ref{st}) hold it
follows that
\begin{equation*}
\left\Vert\Psi(t,t_{0},\psi(t,t_{0},x))\widetilde{R}_{1}(t_{0})v\right\Vert
\leq N_{1}^{x_{0}}e^{\nu_{1}^{x_{0}}}\left\Vert
\widetilde{R}_{1}(t_{0})v\right\Vert
\end{equation*}%
\noindent for all $(t,t_{0})\in\mathcal{T}$ and all $(x,v)\in
\mathcal{Y}$, where
\[
\nu_{1}^{x_{0}}=\frac{\ln c}{t_{0}}, \  N_{1}^{x_{0}}=Ne^{(\omega
+ \nu_{1}^{x_{0}})t_{0}}
\]
and constants $N$ and $\omega$ are given by (\ref{eg}).

From the fact that there exist $t_{0}>0$ and $c\in (0,1)$ such
that inequality (\ref{in}) holds it follows that
\begin{equation*}
\left\Vert \widetilde{R}_{2}(t_{0})v\right\Vert \leq
N_{2}^{x_{0}}e^{-\nu_{2}^{x_{0}}t}\left\Vert
\Psi(t,t_{0},\psi(t,t_{0},x))\widetilde{R}_{2}(t_{0})v\right\Vert
\end{equation*}%
\noindent for all $(t,t_{0})\in\mathcal{T}$ and all $(x,v)\in
\mathcal{Y}$, where
\[
\nu_{2}^{x_{0}}=-\frac{\ln c}{t_{0}}, \
N_{2}^{x_{0}}=\widetilde{N}e^{\widetilde{\omega}t_{0}}.
\]
and constants $\widetilde{N}$ and $\widetilde{\omega}$ are given
by (\ref{ed}).

Hence, the pointwise uniform exponential trichotomy for the
skew-evolution semiflow $\xi$ is proved.
\end{proof}

{\footnotesize

\vspace{5mm}

\noindent\begin{tabular}[t]{ll}

Codru\c{t}a Stoica \\
Institut de Math\' ematiques  \\
Universit\' e Bordeaux 1  \\
France  \\
e-mail: \texttt{codruta.stoica@math.u-bordeaux1.fr}
\end{tabular}}

\end{document}